\documentclass[12pt,reqno]{amsart}

\usepackage{%amssymb, latexsym, amsmath, amscd, array, 
hyperref
}

\usepackage{breakurl}   %this allows for url breaks in pdf

\numberwithin{equation}{section}

\newcommand\R {{\mathbb R}}

\title{A two-track tour of Cauchy's \emph{Cours}}

\author[Mikhail Katz]{Mikhail G. Katz} \address{M.~Katz, Department of
  Mathematics, Bar Ilan University, Ramat Gan 5290002
  Israel}\email{katzmik@math.biu.ac.il}


\begin{document}

\maketitle

\begin{abstract}
Cauchy published his \emph{Cours d'Analyse} 200 years ago.  We analyze
Cauchy's take on the concepts of \emph{rigor}, \emph{continuity}, and
\emph{limit}, and explore a pair of approaches in the literature to
the meaning of his infinitesimal analysis and his sum theorem on the
convergence of series of continuous functions.
\end{abstract}

\section{Rigor then and now}

Building upon pioneering work by Kepler, Fermat, Cavalieri, Gregory,
Wallis, Barrow, and others, Isaac Newton and Gottfried Wilhelm von
Leibniz invented the calculus in the 17th century.  While immediately
acquiring an enthusiastic following, the new methods proved to be
controversial in the eyes of some of their contemporaries, who
employed the more traditional methods of their predecessors.  One of
the controversial aspects of the new technique was Leibniz's
distinction between assignable and inassignable quantities (including
infinitesimals and infinite quantities; see \cite{14c}, \cite{17b},
\cite{21a}).  At the French Academy, the opposition to the new
calculus was led by Michel Rolle, and across the Channel, by George
Berkeley.  The scientific success of the new methods ultimately
silenced the opposition, but lingering doubts persisted (fed in part
by doctrinal theological issues; see \cite{18d}).  A new era was
ushered in by Augustin-Louis Cauchy's textbook \emph{Cours d'Analyse},
addressed to the students of the Ecole Polytechnique in Paris.

Cauchy published his \emph{Cours d'Analyse} (CDA) 200 years ago.  The
book was of fundamental importance for the development of both real
and complex analysis.  Hans Freudenthal mentioned in his essay on
Cauchy (for the \emph{Dictionary of Scientific Biography}) \cite{Fr81}
that Niels Henrik Abel described the CDA as ``an excellent work which
should be read by every analyst who loves mathematical rigor.''  But
what did \emph{rigor} mean to Abel and Cauchy?

In the early 19th century context, the term \emph{rigor} referred to
the standard of mathematical precision set by the geometry of Euclid.
This context enables us to understand Cauchy's intention when, in the
introduction to CDA, he referred to `all the rigor one demands from
geometry', committing himself further to `never rely on arguments
taken from the generality of algebra'.  We see that Cauchy's notion of
rigor in CDA is distinct from ``what has been called the
nineteenth-century `rigorization' of real analysis''
\cite[p.\;221]{Ka21}.

Like the term \emph{rigor}, Cauchy's term \emph{generality of algebra}
requires explanation to be comprehensible to modern readers.  It
refers to certain techniques used by his predecessors, particularly
Euler and Lagrange, that today would be considered cavalier, and
specifically (1)\;proofs based on algebraic manipulation of divergent
series, and (2) the idea that algebraic rules and formulas valid in
the real domain remain valid in the complex domain.  Some of the
proofs that fall under item (1) have since been justified using
summation techniques developed later; some of the techniques under
item (2), in terms of analytic continuation.

By the standards of the current century, some of Freudenthal's
comments could be considered controversial.  Thus, Freudenthal writes:
\begin{enumerate}\item[]
Terms like ``infinitesimally small'' prevail in Cauchy's limit
arguments and epsilontics still looks far away, but there is one
exception.  His proof {\ldots}~of the well-known theorem
\[
\text{ If} \lim_{x\to\infty}(f(x+1)-f(x))=\alpha, \text{ then}
\lim_{x\to\infty}x^{-1} f(x)=\alpha
\]
is a paragon, and the first example, of epsilontics -- the
character~$\varepsilon$ even occurs there.  \cite[p.\;137]{Fr81}
\end{enumerate}
The claim that infinitesimals `prevail' in Cauchyan foundations of
analysis, whereas~$\varepsilon$-$\delta$ arguments `look far away' and
are limited to a small number of exceptions, may surprise a reader
whose perceptions of Cauchyan rigor are influenced by Judith
Grabiner's views \cite{Gr81} and publications that followed them,
especially if they tend to identify \emph{rigor} with the jettisoning
of infinitesimals in favor of~$\varepsilon$-$\delta$ arguments based
on an `algebra of inequalities'.  Some historians today would view
both Freudenthal's and Grabiner's perspectives as outdated.  But the
dual view of Cauchyan analysis persists in the current literature.

\section{Dual tracks}

The track-A view holds that Cauchy, ahead of his time, worked
primarily with an Archi\-me\-dean continuum, and pioneered many of the
techniques that would become known as~$\varepsilon$-$\delta$ in the next
century.

The track-B view holds that Cauchy, like most of his contemporaries
and colleagues at the Ecole Polytechnique, based his analysis
primarily on variable quantities and infinitesimals; see Laugwitz
\cite{La87}, \cite{La89}.

What were then the Cauchyan foundations of analysis?  While we won't
purport to provide a definitive answer in this short note, we will let
Cauchy speak for himself (using the translation by Bradley and
Sandifer \cite{Br09}).  Cauchy writes in the introduction to CDA:
\begin{enumerate}\item[]
In speaking of the continuity of functions, I could not dispense with
a treatment of the principal properties of infinitely small
quantities, properties which serve as the foundation of the
infinitesimal calculus.  \cite[p.\;1]{Br09}
\end{enumerate}
Track-A advocates read this passage as a concession to the management
of the Ecole, and argue that Cauchy ``could not dispense with a
treatment of {\ldots}\;infinitely small quantities'' because of
explicit mandates from the Ecole, against his better judgement.

Track-B advocates read this passage as a recognition by Cauchy (in a
departure from his pre-1820 approaches) that a convincing and
accessible treatment of continuity necessitates infinitesimals, and
note that Cauchy's favorable judgment of infinitesimals is
corroborated by their use in his research long after the end of his
teaching stint at the Ecole; see \cite{20b}.

Readers searching for an~$\varepsilon$-$\delta$ definition of limit in
CDA may be surprised to find instead the following definition:
\begin{enumerate}\item[]
We call a quantity \emph{variable} if it can be considered as able to
take on successively many different values.  {\ldots} When the values
successively attributed to a particular variable indefinitely approach
a fixed value in such a way as to end up by differing from it by as
little as we wish, this fixed value is called the \emph{limit} of all
the other values.  \cite[p.\;6]{Br09}
\end{enumerate}
Here the notion of a variable quantity is taken as primary, and limits
are defined in terms of variable quantities.  Variable quantities
similarly provide the basis for the definition of infinitesimals:
\begin{enumerate}\item[]
When the successive numerical values of such a variable decrease
indefinitely, in such a way as to fall below any given number, this
variable becomes what we call \emph{infinitesimal}, or an
\emph{infinitely small quantity}.  A variable of this kind has zero as
its limit.  \cite[p.\;7]{Br09}
\end{enumerate}
Track-A advocates read this as asserting that an infinitesimal is
merely a null sequence (i.e., a sequence tending to zero), and take
the last sentence to refer to infinitesimals.  Track-B advocates point
out that Cauchy did not write that a variable quantity \emph{is} an
infinitesimal, but rather that a variable quantity \emph{becomes} an
infinitesimal, implying a change in nature (from being a variable
quantity to being an infinitesimal); and take the last sentence to
refer to the variable quantity mentioned at the beginning of the
passage.

Returning to limits, Cauchy writes:
\begin{enumerate}\item[]
When a variable quantity converges towards a fixed limit, it is often
useful to indicate this limit with particular notation. We do this by
placing the abbreviation
\[
\text{lim}
\]
in front of the variable quantity in question.  \cite[p.\;12]{Br09}
\end{enumerate}
Bradley and Sandifer note that the 1821 edition of the CDA used the
notation ``lim.''\;(with a period).  Cauchy's description of
``lim.''~as an abbreviation suggests that he viewed it as secondary
to, or an aspect of, the concept of a variable quantity.

In Chapter 2, Cauchy returns to the definition of infinitesimals:
\begin{enumerate}\item[]
We say that a variable quantity becomes \emph{infinitely small} when
its numerical value decreases indefinitely in such a way as to
converge towards the limit zero.  \cite[p.\;21]{Br09}
\end{enumerate}
The relation between the concepts of variable quantity and
infinitesimal was already discussed above, as well as the possible
ambiguity of the verb \emph{becomes}; limits again play a secondary
role.  Cauchy proceeds next to the properties of infinitesimals:
\begin{enumerate}\item[]
Infinitely small and infinitely large quantities enjoy several
properties that lead to the solution of important questions, which I
will explain in a few words.  Let~$\alpha$ be an infinitely small
quantity, that is a variable whose numerical value decreases
indefinitely.  \cite[p.\;22]{Br09}
\end{enumerate}
Track-A advocates point out that here Cauchy states that an
infinitesimal is a variable quantity.  Track-B advocates note that
here Cauchy is no longer dealing with detailed definitions, and this
particular formulation is merely shorthand for the more careful
definition in terms of \emph{becoming} elaborated earlier.

\section{Continuity in 1817, 1821, and beyond}

Cauchy's first documented characterisation of continuity is found in a
record of a course summary dating from march 1817 (a month before the
earliest written mention of Bolzano's \emph{Rein analytischer Beweis}
in an Olms catalogue).  The definition can be described as reasonably
precise in the sense of enabling a straightforward transcription as an
impeccable modern definition; see~\cite{20a}.  In modern mathematics,
a real function~$\varphi$ is continuous at~$c\in\R$ if and only if for
each sequence~$(x_n)$ converging to~$c$, one has
\begin{equation}
\label{e01}
\varphi(c)= \varphi\Big(\lim_{n\to\infty}x_n\Big)=
\lim_{n\to\infty}\varphi(x_n),
\end{equation}
or briefly~$\varphi\circ \lim=\lim\circ\, \varphi$ at~$c$, expressing
the commutation of $\varphi$ and $\lim$.  In 1817, Cauchy wrote:
\begin{enumerate}\item[]
The limit of a continuous function of several variables is [equal to]
the same function of their limit.\, Consequences of this Theorem with
regard to the continuity of composite functions dependent on a single
variable.
\end{enumerate}
(Being part of a summary, the second phrase is not a complete
sentence.)  Cauchy's 1817 characterisation of continuity in terms of
commutation of~$\varphi$ and~$\lim$ as in \eqref{e01} does not use
infinitesimals and thus contrasts with his definitions involving
infinitesimals given four years later in CDA.\, Surprisingly, it is
the 1817 characterisation that is actually used in CDA; see
Section~\ref{s04}.

Here is Cauchy's first definition of continuity in CDA:
\begin{enumerate}\item[]
Let~$f(x)$ be a function of the variable~$x$, and suppose that for
each value of~$x$ between two given limits, the function always takes
a unique finite value.  If, beginning with a value of~$x$ contained
between these limits, we add to the variable~$x$ an \emph{infinitely
  small increment}~$\alpha$, the function itself is incremented by the
difference \hbox{$f(x+\alpha)-f(x)$,} which depends both on the new
variable~$\alpha$ and on the value of~$x$.  Given this, the
function~$f(x)$ is a \emph{continuous} function of~$x$ between the
assigned limits if, for each value of~$x$ between these limits, the
numerical value of the difference~$f(x+\alpha)-f(x)$ decreases
indefinitely with the numerical value of~$\alpha$.
\cite[p.\;26]{Br09}
\end{enumerate}
Note that, while the increment~$\alpha$ is described as infinitesimal,
the resulting change~$f(x+\alpha)-f(x)$ is not.  This 1821 definition
can be seen as intermediary between the march\;1817 characterisation
in terms of variables (not mentioning infinitesimals), and his second
1821 definition stated purely in terms of infinitesimals:
\begin{enumerate}\item[]
In other words, \emph{the function~$f(x)$ is continuous with respect
  to~$x$ between the given limits if, between these limits, an
  infinitely small increment in the variable always produces an
  infinitely small increment in the function itself}.  (ibid.)
\end{enumerate}
Significantly, Cauchy's abbreviation ``lim.''~appeared in none of the
definitions of continuity given in CDA.%
\footnote{Surprisingly, a contrary claim appears in the recent
  literature: ``Cauchy gave a faultless definition of continuous
  function, using the notion of `limit' for the first time.  Following
  Cauchy's idea, Weierstrass popularized the $\epsilon$-$\delta$
  argument in the 1870's'' Dani--Papadopoulos \cite{Da19}, 2019,
  p.\;283).  The notion of limit does not appear in any of Cauchy's
  definitions of continuity in CDA.}
It is the second definition purely in terms of infinitesimals that
reappears in the following works by Cauchy:
\begin{itemize}
\item
\emph{R\'esum\'e des Le\c{c}ons} (1823); English translation
\cite[p.\;9]{Ca19};
\item
\emph{Le\c{c}ons sur le Calcul Diff\'erentiel} (1829)
\cite[p.\;9]{Ca29};
\item
\emph{M\'emoire sur l'Analyse Infinit\'esimale} (1844)
\cite[p.\;17]{Ca44};
\item
the 1853 article \cite{Ca53} on the sum theorem (see
Section~\ref{s4}).
\end{itemize}
Nonetheless, the perception that Cauchy allegedly `sought to establish
foundations for real analysis that gave no role to infinitesimals' has
firmly entered the canonical creation narrative of modern mathematical
analysis; see e.g., \cite{Op21}.

\section{A modern digression: infinitesimals without Choice}
\label{s03}

Cauchy's definition of continuity in CDA, in its B-track
interpretation, is harder to follow for modern readers more familiar
with the~$\varepsilon$-$\delta$ definition of continuity \`a la
Weierstrass--Dini than with the definition in a modern infinitesimal
theory.  We therefore provide a formalisation of Cauchy's procedures
involving continuity in terms of the theory SPOT (acronym of its
axioms) developed in \cite{21d}; see also \cite{17f}.  SPOT has the
advantage of being conservative over the traditional Zermelo--Fraenkel
set theory (ZF) and therefore depends on neither the Axiom of Choice
nor the existence of ultrafilters.

The language of ZF is limited to the two-place membership
relation~$\in$.  The language of SPOT includes also a predicate ST,
where ST$(x)$ reads `$x$ is standard'.  Such a distinction between
standard and nonstandard entities can be thought of as formalizing the
Leibnizian distinction between assignable and inassignable quantities;
see \cite{21a}.  The standard ordered field~$\R$ has both standard and
nonstandard elements.  An element~$\alpha$ is \emph{infinitesimal}
if~$|\alpha|<r$ for each standard~$r>0$.  Let $x$ be a standard point
in the domain of a real standard function $f$.  Then~$f$ is continuous
at~$x$ (in the traditional sense of the~$\varepsilon$-$\delta$
definition) if and only if
\begin{equation}
\label{e02}
\text{infinitesimal } \alpha \text{ produce infinitesimal changes }
f(x+\alpha)-f(x),
\end{equation}
whenever $x+\alpha$ is in the domain of $f$.  Continuity in an
interval, say~$(0,1)$, is equivalent to the satisfaction of
condition~\eqref{e02} at every standard point~$x\in(0,1)$.  

It should be mentioned that, according to most scholars, Cauchy did
not formulate the notion of continuity at a point, but only
`continuity between limits' (i.e., in an interval) or in a
neighborhood of a point.  Freudenthal remarks:
\begin{enumerate}\item[]
It is the weakest point in Cauchy's reform of calculus that he never
grasped the importance of uniform continuity.  \cite[p.\;137]{Fr81}
\end{enumerate}

From an A-track viewpoint, the weakness is that there seems to be no
trace in Cauchy of the idea that a significant issue is whether the
allowable error is independent of the point~$x$ or not.

From a B-track viewpoint, the weakness is that Cauchy did not make it
clear whether condition~\eqref{e02} is expected to be satisfied only
at assignable points~$x$ or at all points of the interval.  Note that
uniform continuity of $f$ on, say, $(0,1)$ is equivalent to
\eqref{e02} being satisfied at all points of~$(0,1)$.  For example,
$\frac1x$ fails to be uniformly continuous because of the failure of
\eqref{e02} at an infinitesimal input $x=\beta>0$: indeed, the change
$\frac1{\beta+\alpha}-\frac1{\beta}$ is not infinitesimal if, say,
$\alpha=\beta$.  See further in Section~\ref{s4}.

By Section 2.3, Cauchy reaches the Theorem described by Freudenthal as
a \emph{paragon} of~$\varepsilon$-$\delta$ arguments.  What Cauchy
actually shows is that if~$f(x+1)-f(x)$ is between~$k-\varepsilon$ and
$k+\varepsilon$ then (assuming monotonicity)~$\frac{f(x)}{x}$ is
similarly between~$k-\varepsilon$ and $k+\varepsilon$.  If anything
this is a paragon of~$\varepsilon$-$\varepsilon$ arguments, since here
the~$\delta$ equals $\varepsilon$!  The trademark feature of
modern~$\varepsilon$-$\delta$ arguments, namely an explicit
(nontrivial) dependence of~$\delta$ on~$\varepsilon$, does not appear
anywhere in Cauchy's alleged
%
%e.g. in Sutherland
%
`algebra of inequalities', lending support to Freudenthal's sentiment
that `epsilontics looks far away'.

Here Cauchy mentions that an infinite limit `is larger than any
\emph{assignable} number' indicating familiarity with this Leibnizian
term (which occurs nine times in CDA).

\section{Functional equations and continuity via commutation}
\label{s04}

In Chapter 5 Cauchy studies functional relations for continuous
functions, and treats the following problem:
\begin{enumerate}\item[]
Problem I. -- To determine the function~$\varphi (x)$ in such a manner
that it remains continuous between any two real limits of the variable
$x$ and so that for all real values of the variables~$x$ and~$y$, we
have~$\varphi (x+y) = \varphi (x)+\varphi (y)$.  \cite[p.\;71]{Br09}
\end{enumerate}
Cauchy arrives at
$\varphi\big(\frac{m}{n}\alpha\big)=\frac{m}{n}\varphi(\alpha)$ and
argues as follows:
\begin{enumerate}\item[]
Then, by supposing that the fraction~$\frac{m}{n}$ varies in such a
way as to converge towards any number~$\mu$, and passing to the limit,
we find that~$\varphi (\mu\alpha) = \mu\varphi(\alpha)$.
\cite[p.\;72]{Br09}
\end{enumerate}
Here Cauchy exploits the 1817 characterisation of the continuity
of~$\varphi$ in terms of the commutation of~$\varphi$ and~$\lim$ as
summarized in formula~\eqref{e01}, rather than the definitions
presented in Chapter 2 of CDA (the 1817 characterisation is also used
in the proof of the Intermediate Value Theorem).  Curiously, Cauchy
provides no explanatory comment.  Possibly, Cauchy wrote the material
in Chapter\;5 with definition~\eqref{e01} in mind, and introduced the
definitions in Chapter 2 at a later stage in the writing of the book.
Cauchy applies a similar technique to study the functional
relation~$\varphi(x+y)=\varphi(x)\varphi(y)$ and other variations.

\section
{Sum theorem and convergence \emph{always}}
\label{s4}

Chapter 6 includes Cauchy's controversial sum theorem.  We summarize
the historical facts.  The 1821 formulation of the theorem appears to
be incorrect to the modern reader, as it seems to assert that
pointwise convergence of a series of continuous functions implies the
continuity of the sum.  Already in 1826 Abel pointed out that the
theorem `suffers exceptions'.  Cauchy was curiously silent on the
subject of the sum theorem for several decades.  Then in 1853 he
presented a modified statement of the sum theorem, mentioned an
example similar to Abel's (without mentioning Abel by name), and
explained why the example does not contradict the (modified) theorem.
Numerous scholars have attempted to explain what the modification was
(if any) and to interpret Cauchy's sum theorem in terms intelligible
to modern audiences.

Cauchy considers the series obtained as the sum of the terms of the
sequence $u_0,u_1,u_2,\ldots,u_n,u_{n+1},\ldots$, denoted (1).  In its
1821 formulation, the theorem asserts the following.

\begin{enumerate}\item[]
When the various terms of series (1) are functions of the same
variable~$x$, continuous with respect to this variable in the
neighborhood of a particular value for which the series converges, the
sum~$s$ of the series is also a continuous function of~$x$ in the
neighborhood of this particular value.  \cite[p.\;90]{Br09}
\end{enumerate}
Here is the 1853 formulation:

\begin{enumerate}\item[]
%
%Si les diff\'erents termes de la s\'erie
%
If the various terms of the series
\[
\tag{1}
\hfill u_0, u_1, u_2, \ldots, u_n, u_{n+1}, \ldots \hfil
\]
%
%sont des fonctions de la variable r\'eelle~$x$, continues, par rapport
%\`a cette variable, entre des limites donn\'ees; si, d'ailleurs, la
%somme
%
are functions of a real variable $x$ which are continuous, with
respect to this variable, between the given bounds; and if,
furthermore, the sum
\begin{equation*}
\tag{3}
u_n+u_{n+1} + \ldots + u_{n'-1}
\end{equation*}
%
%devient \emph{toujours} infiniment petite pour des valeurs infiniment
%grandes des nombres entiers~$n$ et~$n'>n$, la s\'erie~(1) sera
%convergente et la somme~$s$ de la s\'erie sera, entre les limites
%donn\'ees, fonction continue de la variable~$x$.  
%
\emph{always} becomes infinitely small for infinitely large values of
the whole numbers $n$ and $n'>n$, then the series (1) will converge
and the sum $s$ of the series will be, between the given bounds, a
continuous function of the variable $x$.  \cite[pp.\;456--457]{Ca53}
(translation ours)
\end{enumerate}
Note that Cauchy adds the word \emph{toujours} (always).  But what
exactly is supposed to happen always in 1853, and how does this modify
the 1821 hypothesis?  Cauchy himself provides a hint in his 1853
analysis of the example $\sum_n \frac{\sin nx}{n}$, representing a
(discontinuous) sawtooth waveform.  But the hint he provides is itself
puzzling.  Cauchy evaluates $u_n+u_{n+1} + \ldots+ u_{n'-1}$ at
$x=\frac1n$ and shows that the sum does not become arbitrarily small,
and in fact can be made ``sensibly equal'' to the integral
$\int_1^\infty \frac{\sin x}{x}dx=0.6244\ldots$

Here Cauchy explicitly describes $n$ (and $n'$) as infinite; then
$x=\frac1n$ is infinitesimal.

Many scholars of both A-track and B-track persuasion have argued that
Cauchy meant to add what is known today as the condition of uniform
convergence.

A-track advocates may interpret the condition as requiring the sums
$u_n+u_{n+1} + \ldots + u_{n'-1}$ to be small independently of~$n,n'$
and also of the input~$x$, a condition that can be stated formally in
terms of an alternating quantifier string of the type
\[
\forall\varepsilon>0\;\exists\delta>0\;\forall
n<n'\in\mathbb{N}\;\forall x\;\big(n>\tfrac1\delta\to
|u_n(x)+u_{n+1}(x)+\ldots+u_{n'-1}(x)|<\varepsilon\big).
\]
Something along these (long) lines would have to be attributed to
Cauchy, in inchoate form, in order to interpret the addition of
uniform convergence in an A-track fashion.  What is unclear is how the
word \emph{always} manages to allude to such independence,
particularly since Cauchy seems to have overlooked its significance in
the context of continuity (see Freudenthal's remark quoted in
Section~\ref{s03} and the ensuing discussion).  Moreover, the
evaluation at what seems to be a new type of point, namely
$x=\frac1n$, appearing in Cauchy's discussion of $\sum_n \frac{\sin
  nx}{n}$, suggests that the insistence on the qualifier \emph{always}
indicates an extension of the inputs $x$ to include additional ones
(that were not \emph{always} included before).

A B-track reading of Cauchy's 1853 hypothesis interprets the qualifier
\emph{always} as referring to additional inassignable inputs $x$
(including infinitesimal values such as $\frac1n$ for infinite $n$).
If one requires the sum $u_n+u_{n+1} + \ldots+ u_{n'-1}$ to be
infinitesimal for all infinite $n,n'$ and all inputs $x$ (standard and
nonstandard) then one indeed obtains a condition equivalent to uniform
convergence (Robinson \cite[Theorem\;4.6.1, p.\;116]{Ro66}),
guaranteeing the continuity of the limit function.  For more details
see \cite{18e}, \cite{19a}.

\section{Conclusion}

Though we have sketched widely divergent readings of Cauchy's
definitions and his theorems, we hope to have conveyed to the reader a
sense that not only the concept of \emph{rigor} but also
\emph{continuity} and \emph{limit} may have had a different meaning to
Cauchy than they do to us today, underscoring the contingency of the
historical evolution of mathematics.  While Cauchy incontestably made
extensive use of \emph{bona fide} infinitesimals in fields as varied
as differential geometry, elasticity theory, and geometric probability
(see \cite{20b}), the nature of his foundational stance remains
controversial.  But perhaps this is as it should be: with a genius of
Cauchy's caliber, tidy construals of his work may necessarily amount
to a flattening of his multi-dimensional vision.

\end{document}